\documentclass[leqno,draft]{article}

\newtheorem{theorem}{Theorem}
\newtheorem{lemma}[theorem]{Lemma}
\newtheorem{proposition}[theorem]{Proposition}
\newtheorem{definition}[theorem]{Definition}
\newtheorem{corollary}[theorem]{Corollary}

\newcommand{\begintheorem}{\addtocounter{equation}{1}\begin{theorem}}
\newcommand{\beginlemma}{\addtocounter{equation}{1}\begin{lemma}}
\newcommand{\beginproposition}{\addtocounter{equation}{1}\begin{proposition}}
\newcommand{\begindefinition}{\addtocounter{equation}{1}\begin{definition}}
\newcommand{\begincorollary}{\addtocounter{equation}{1}\begin{corollary}}

\begin{document}

\title{Some remarks about solenoids, 2}

\author{Stephen Semmes \\
        Rice University}

\date{}

\maketitle

\begin{abstract}
A class of solenoids is considered, including some aspects in $n$
(topological) dimensions, where one basically gets some fractal
versions of tori.
\end{abstract}

\tableofcontents

\section{$p$-Adic numbers}
\label{p-adic numbers}
\setcounter{equation}{0}

        If $p$ is a prime number, then the \emph{$p$-adic absolute value}
$|x|_p$ is defined on the field ${\bf Q}$ of rational numbers as follows.
Of course, $|x|_p = 0$ when $x = 0$.  Otherwise, $x$ can be expressed
as $p^j \, (a/b)$ for some integers $a, b, j$ such that $a, b \ne 0$
and neither $a$ nor $b$ is divisible by $p$.  In this case,
\begin{equation}
\label{|x|_p = p^{-j}}
        |x|_p = p^{-j}.
\end{equation}
It is well known and easy to check that
\begin{equation}
\label{|x + y|_p le max(|x|_p, |y|_p)}
        |x + y|_p \le \max(|x|_p, |y|_p)
\end{equation}
and
\begin{equation}
\label{|x y|_p = |x|_p |y|_p}
        |x \, y|_p = |x|_p \, |y|_p
\end{equation}
for every $x, y \in {\bf Q}$.  The \emph{$p$-adic metric} on ${\bf Q}$
is defined by
\begin{equation}
\label{d_p(x, y) = |x - y|_p}
        d_p(x, y) = |x - y|_p.
\end{equation}
This is actually an ultrametric on ${\bf Q}$, since it satisfies
the ultrametric version of the triangle inequality
\begin{equation}
\label{d_p(x, z) le max(d_p(x, y), d_p(y, z))}
        d_p(x, z) \le \max(d_p(x, y), d_p(y, z))
\end{equation}
for every $x, y, z \in {\bf Q}$, by (\ref{|x + y|_p le max(|x|_p, |y|_p)}).

        The field ${\bf Q}_p$ of \emph{$p$-adic numbers} is obtained by
completing ${\bf Q}$ with respect to the $p$-adic metric, in the same
way that the field ${\bf R}$ of real numbers may be obtained by
completing ${\bf Q}$ with respect to the standard metric $|x - y|$,
associated to the ordinary absolute value $|x|$.  In particular, the
$p$-adic absolute value $|x|_p$ has a natural extension to $x \in {\bf
  Q}_p$, such that $|x|_p$ is still an integer power of $p$ when $x
\ne 0$, and (\ref{|x + y|_p le max(|x|_p, |y|_p)}) and (\ref{|x y|_p =
  |x|_p |y|_p}) hold for every $x, y \in {\bf Q}_p$.  The natural extension
of the $p$-adic metric to ${\bf Q}_p$ corresponds to the extension of the
$p$-adic absolute value as in (\ref{d_p(x, y) = |x - y|_p}), and is an
ultrametric on ${\bf Q}_p$.  Of course, ${\bf Q}$ is dense in ${\bf Q}_p$
with respect to the $p$-adic metric, by construction.

        The set ${\bf Z}_p$ of $p$-adic integers may be defined by
\begin{equation}
\label{{bf Z}_p = {x in {bf Q}_p : |x|_p le 1}}
        {\bf Z}_p = \{x \in {\bf Q}_p : |x|_p \le 1\}.
\end{equation}
Equivalently, ${\bf Z}_p$ is the closed unit ball in ${\bf Q}_p$, and
hence a closed set with respect to the $p$-adic metric.  It is easy to
see that ${\bf Z}_p$ is a sub-ring of ${\bf Q}_p$, using (\ref{|x +
  y|_p le max(|x|_p, |y|_p)}) and (\ref{|x y|_p = |x|_p |y|_p}).
Clearly ${\bf Z}_p$ contains the usual ring ${\bf Z}$ of integers,
since $|x|_p \le 1$ when $x \in {\bf Z}$ by the definition of $|x|_p$.
One can check that ${\bf Z}$ is dense in ${\bf Z}_p$ with respect to
the $p$-adic metric, so that ${\bf Z}_p$ is the same as the closure
of ${\bf Z}$ in ${\bf Q}_p$.  More precisely, one can first check that
${\bf Z}$ is dense in
\begin{equation}
\label{{bf Z}_p cap {bf Q} = {x in {bf Q} : |x|_p le 1}}
        {\bf Z}_p \cap {\bf Q} = \{x \in {\bf Q} : |x|_p \le 1\}
\end{equation}
with respect to the $p$-adic metric, so that (\ref{{bf Z}_p cap {bf Q}
  = {x in {bf Q} : |x|_p le 1}}) is the same as the closure of ${\bf
  Z}$ in ${\bf Q}$ with respect to $d_p(x, y)$.  Every $x \in {\bf
  Z}_p$ can be approximated by $y \in {\bf Q}$ with respect to the
$p$-adic metric, because ${\bf Q}$ is dense in ${\bf Q}_p$, and $|y|_p
\le 1$ when $|x - y|_p \le 1$, by (\ref{|x + y|_p le max(|x|_p,
  |y|_p)}).  Thus (\ref{{bf Z}_p cap {bf Q} = {x in {bf Q} : |x|_p le
    1}}) is dense in ${\bf Z}_p$, which implies that ${\bf Z}$ is
dense in ${\bf Z}_p$, since ${\bf Z}$ is dense in (\ref{{bf Z}_p cap
  {bf Q} = {x in {bf Q} : |x|_p le 1}}).

        Similarly,
\begin{equation}
\label{p^l {bf Z}_p = ... = {y in {bf Q}_p : |y|_p le p^{-l}}}
        p^l \, {\bf Z}_p = \{p^l \, x : x \in {\bf Z}_p\}
                         = \{y \in {\bf Q}_p : |y|_p \le p^{-l}\}
\end{equation}
is a closed sub-ring of ${\bf Q}_p$ for each integer $l$, and an ideal
in ${\bf Z}_p$ when $l \ge 0$.  Of course, $p^l \, {\bf Z}$ is also an
ideal in ${\bf Z}$ when $l \ge 0$, and there is a natural ring
homomorphism from ${\bf Z} / p^l \, {\bf Z}$ into ${\bf Z}_p / p^l \,
{\bf Z}_p$, corresponding to the inclusion of ${\bf Z}$ in ${\bf
  Z}_p$.  It is easy to see that this homomorphism is injective,
because
\begin{equation}
\label{{bf Z} cap p^l {bf Z}_p = p^l {bf Z}}
        {\bf Z} \cap p^l \, {\bf Z}_p = p^l \, {\bf Z}
\end{equation}
for each nonnegative integer $l$.  One can also check that this
homomorphism is surjective, using the fact that ${\bf Z}$ is dense in
${\bf Z}_p$.  Thus ${\bf Z}_p / p^l \, {\bf Z}_p$ is isomorphic as a
ring to ${\bf Z} / p^l \, {\bf Z}$ for each nonnegative integer $l$.

        In particular, ${\bf Z}_p / p^l \, {\bf Z}_p$ has exactly
$p^l$ elements for each nonnegative integer $l$.  Equivalently,
${\bf Z}_p$ can be expressed as the union of $p^l$ pairwise-disjoint
translates of $p^l \, {\bf Z}_p$ for each $l \ge 0$.  This implies
that ${\bf Z}_p$ is totally bounded in ${\bf Q}_p$, which means that
it can be covered by finitely many balls of arbitrarily small radius.
It is well known that a subset of a complete metric space is compact
if and only if it is closed and totally bounded, so that ${\bf Z}_p$ 
is compact with respect to the $p$-adic metric.  Similarly, 
$p^j \, {\bf Z}_p$ is a compact set in ${\bf Q}_p$ for each integer $j$.

        By standard arguments, addition and multiplication on ${\bf Q}_p$
define continuous mappings from ${\bf Q}_p \times {\bf Q}_p$ into ${\bf Q}_p$,
and $x \mapsto 1/x$ is continuous as a mapping from 
${\bf Q}_p \backslash \{0\}$ into itself.  Thus ${\bf Q}_p$ is a topological
field with respect to the topology determined by the $p$-adic metric,
and a topological group with respect to addition in particular.
Of course, ${\bf Q}_p$ is locally compact, as in the previous paragraph,
and hence is also equipped with Haar measure $H_p$.  More precisely,
$H_p$ is a nonnegative Borel measure on ${\bf Q}_p$ which is invariant
under translations and satisfies $H_p(U) > 0$ for every nonempty open
set $U \subseteq {\bf Q}_p$ and $H_p(K) < \infty$ when $K \subseteq {\bf Q}_p$
is compact.  Haar measure is unique up to multiplication by a positive
real number, and in this case it is convenient to choose $H_p$ so that
\begin{equation}
        H_p({\bf Z}_p) = 1.
\end{equation}
As usual, Haar measure satisfies some additional regularity
properties, which are automatic in this situation, because open
subsets of ${\bf Q}_p$ are $\sigma$-compact.  This uses the fact that
${\bf Q}_p$ is $\sigma$-compact, as in the previous paragraph, and that
open subsets of any metric space are $F_\sigma$ sets.

        If $j$ is a positive integer, then it follows that
\begin{equation}
\label{H_p(p^j {bf Z}_p) = p^{-j}}
        H_p(p^j \, {\bf Z}_p) = p^{-j},
\end{equation}
because ${\bf Z}_p$ is the union of $p^j$ pairwise-disjoint translates
of $p^j \, {\bf Z}_p$.  The same relation holds when $j$ is a negative
integer, since $p^j \, {\bf Z}_p$ is then the union of $p^{-j}$
pairwise-disjoint translates of ${\bf Z}_p$.  If $E \subseteq {\bf
  Q}_p$ is a Borel set and $a \in {\bf Q}_p$, then $a \, E = \{a \, x
: x \in E\}$ is a Borel set in ${\bf Q}_p$, and
\begin{equation}
\label{H_p(a E) = |a|_p H_p(E)}
        H_p(a \, E) = |a|_p \, H_p(E).
\end{equation}
Note that $H_p$ is the same as one-dimensional Hausdorff measure on
${\bf Q}_p$ with respect to the $p$-adic metric.

\section{Finitely many primes}
\label{finitely many primes}
\setcounter{equation}{0}

        Let $S = \{p_1, \ldots, p_k\}$ be a nonempty finite set of
prime numbers, and let ${\bf Z}_S = {\bf Z}[1/p_1, \ldots, 1/p_k]$
be the ring generated by their multiplicative inverses.  Thus
the elements of ${\bf Z}_S$ are rational numbers of the form $a / b$,
where $a$ and $b$ are integers, and $b \ne 0$ is a product of powers
of elements of $S$.  Let $R_S$ be the Cartesian product
\begin{equation}
        R_S = {\bf Q}_{p_1} \times \cdots \times {\bf Q}_{p_k} \times {\bf R}
\end{equation}
equipped with the product topology, using the topology on ${\bf
  Q}_{p_j}$ associated to the $p_j$-adic metric for each $j = 1,
\ldots, k$, and the standard topology on the real line ${\bf R}$.
Of course, $R_S$ is also a commutative ring, where addition and
multiplication are defined coordinatewise, and more precisely a 
topological ring with respect to the topology just mentioned.

        Because ${\bf Q} \subseteq {\bf Q}_{p_j}$ for $j = 1, \ldots, k$
and ${\bf Q} \subseteq {\bf R}$, there is a natural diagonal embedding
$\Delta_S$ of ${\bf Q}$ in $R_S$, which sends each $x \in {\bf Q}$ to
the element of $R_S$ whose projections in the ${\bf Q}_{p_j}$'s and ${\bf R}$
are all equal to $x$.  Note that $\Delta_S$ is a ring homomorphism
from ${\bf Q}$ into $R_S$, so that $\Delta_S({\bf Q})$ is a sub-ring
of $R_S$.  Similarly, ${\bf Z}_S$ is a sub-ring of ${\bf Q}$, and hence
$\Delta_S({\bf Z}_S)$ is also a sub-ring of $R_S$.  It is well known that
$\Delta_S({\bf Z}_S)$ is discrete in $R_S$, in much the same way that
${\bf Z}$ is discrete in ${\bf R}$.  To see this, observe first that
if $x \in {\bf Z}_S$ satisfies
\begin{equation}
\label{|x|_{p_j} le 1}
        |x|_{p_j} \le 1
\end{equation}
for each $j = 1, \ldots, k$, then $x \in {\bf Z}$.  If we also have that
\begin{equation}
\label{|x| < 1}
        |x| < 1,
\end{equation}
then it follows that $x = 0$, so that $0$ is the only element of
$\Delta_S({\bf Z}_S)$ in a suitable neighborhood of $0$ in $R_S$.
Similarly, if $y \in {\bf Z}_S$, then there is a neighborhood of
$\Delta_S(y)$ in $R_S$ such that $y$ is the only element of
$\Delta_S({\bf Z}_S)$ in that neighborhood, by translation-invariance.

        Alternatively, let us denote elements of $R_S$ by 
$r = (r_1, \ldots, r_k, r_{k + 1})$, where $r_j \in {\bf Q}_{p_j}$
for $j = 1, \ldots, k$, and $r_{k + 1} \in {\bf R}$.  If $r, r' \in R_S$,
then put
\begin{equation}
\label{d_S(r, r') = max(|r_1 - r'_1|_{p_1}, ..., |r_{k + 1} - r'_{k + 1}|)}
        d_S(r, r') = \max(|r_1 - r'_1|_{p_1}, \ldots, |r_k - r'_k|_{p_k}, 
                                                     |r_{k + 1} - r'_{k + 1}|).
\end{equation}
This defines a metric on $R_S$, which determines the same topology on
$R_S$ as before, and which is invariant under translations on $R_S$.
The argument in the previous paragraph implies that
\begin{equation}
\label{d_S(Delta_S(x), Delta_S(y)) ge 1}
        d_S(\Delta_S(x), \Delta_S(y)) \ge 1
\end{equation}
for every $x, y \in {\bf Z}_S$ with $x \ne y$.  Let us check that each
element of $R_S$ is at a bounded distance from an element of
$\Delta_S({\bf Z}_S)$.

        Let $r \in R_S$ be given.  Suppose first that $1 \le j \le k$, 
so that $r_j \in {\bf Q}_{p_j}$.  If $|r_j|_{p_j} \le 1$, then put
$x_j = 0$.  Otherwise, $|r_j|_{p_j} = p_j^{l_j}$ for some positive
integer $l_j$, which implies that
\begin{equation}
\label{|p_j^{l_j} r_j|_{p_j} = 1}
        |p_j^{l_j} \, r_j|_{p_j} = 1.
\end{equation}
Because ${\bf Z}$ is dense in ${\bf Z}_{p_j}$, there is an integer
$a_j$ such that
\begin{equation}
\label{|p_j^{l_j} r_j - a_j|_{p_j} le p_j^{-l_j}}
        |p_j^{l_j} \, r_j - a_j|_{p_j} \le p_j^{-l_j},
\end{equation}
or equivalently,
\begin{equation}
\label{|r_j - a_j p^{-l_j}|_{p_j} le 1}
        |r_j - a_j \, p^{-l_j}|_{p_j} \le 1.
\end{equation}
In this case, we take $x_j = a_j \, p_j^{-l_j}$, which also satisfies
\begin{equation}
\label{|x_j|_p le 1}
        |x_j|_p \le 1
\end{equation}
for each prime $p \ne p_j$.  Thus $x_j \in {\bf Z}_S$ in both cases,
and hence
\begin{equation}
\label{x = sum_{j = 1}^k x_j}
        x = \sum_{j = 1}^k x_j
\end{equation}
is an element of ${\bf Z}_S$ too.  It is easy to see that
\begin{equation}
\label{|r_j - x|_{p_j} le 1}
        |r_j - x|_{p_j} \le 1,
\end{equation}
when $1 \le j \le k$, because $|r_j - x_j|_{p_j} \le 1$, and
$|x_i|_{p_j} \le 1$ when $i \ne j$.  If $y \in {\bf Z}_S$ has the
property that $y - x \in {\bf Z}$, then
\begin{equation}
\label{|r_j - y|_{p_j} le 1}
        |r_j - y|_{p_j} \le 1
\end{equation}
for $j = 1, \ldots, k$ as well.  This permits us to choose $y$ such
that $y - x \in {\bf Z}$ and
\begin{equation}
\label{|r_{k + 1} - y| le 1/2}
        |r_{k + 1} - y| \le 1/2,
\end{equation}
since $r_{k + 1} \in {\bf R}$.  It follows that
\begin{equation}
\label{d_S(r, Delta_S(y)) le 1}
        d_S(r, \Delta_S(y)) \le 1,
\end{equation}
as desired.

\section{Snowflake metrics}
\label{snowflake metrics}
\setcounter{equation}{0}

        Let $p$ be a prime number, and let $a$ be a positive real number.
If we put
\begin{equation}
\label{|x|_{p, a} = |x|_p^a}
        |x|_{p, a} = |x|_p^a
\end{equation}
for each $x \in {\bf Q}_p$, then we get another ``absolute value function''
on ${\bf Q}_p$, which still satisfies
\begin{equation}
\label{|x + y|_{p, a} le max(|x|_{p, a}, |y|_{p, a})}
        |x + y|_{p, a} \le \max(|x|_{p, a}, |y|_{p, a})
\end{equation}
and
\begin{equation}
\label{|x y|_{p, a} = |x|_{p, a} |y|_{p, a}}
        |x \, y|_{p, a} = |x|_{p, a} \, |y|_{p, a}
\end{equation}
for every $x, y \in {\bf Q}_p$.  In particular,
\begin{equation}
\label{d_{p, a}(x, y) = |x - y|_{p, a} = |x - y|_p^a = d_p(x, y)^a}
        d_{p, a}(x, y) = |x - y|_{p, a} = |x - y|_p^a = d_p(x, y)^a
\end{equation}
is a translation-invariant ultrametric on ${\bf Q}_p$ that determines
the same topology on ${\bf Q}_p$.  It is easy to see that
$t$-dimensional Hausdorff measure on ${\bf Q}_p$ with respect to
$d_{p, a}(x, y)$ is the same as $(t/a)$-dimensional Hausdorff measure
on ${\bf Q}_p$ with respect to $d_p(x, y)$ for each $t \ge 0$.  Thus
normalized Haar measure on ${\bf Q}_p$ is the same as
$(1/a)$-dimensional Hausdorff measure on ${\bf Q}_p$ with respect to
$d_{p, a}(x, y)$, so that ${\bf Q}_p$ has Hausdorff dimension $1/a$
with respect to $d_{p, a}(x, y)$.

        Similarly, one can check that $|x - y|^a$ defines a metric on the 
real line when $0 < a \le 1$, which determines the same topology on
${\bf R}$.  Although $|x - y|^a$ is not a metric on ${\bf R}$ when $a
> 1$, it is a ``quasi-metric'', which is sufficient for some purposes.
As before, $t$-dimensional Hausdorff measure on ${\bf R}$ with respect
to $|x - y|^a$ is the same as $(t/a)$-dimensional Hausdorff measure on
${\bf R}$ with respect to $|x - y|$.  It is well known that the topological
dimension of a metric space is always less than or equal to the Hausdorff
dimension, and of course the topological dimension of the real line is equal
to $1$.  The geometry of ${\bf R}$ with respect to $|x - y|^a$ is like that
of a snowflake curve when $a < 1$, with Hausdorff dimension $1/a$.

        Let us return now to the situation of the previous section,
and let $a = (a_1, \ldots, a_k)$ be a $k$-tuple of real numbers.  As
an extension of (\ref{d_S(r, r') = max(|r_1 - r'_1|_{p_1}, ..., |r_{k
    + 1} - r'_{k + 1}|)}), put
\begin{equation}
\label{d_{S, a}(r, r') = max(...)}
 d_{S, a}(r, r') = \max(|r_1 - r_1'|_{p_1}^{a_1}, \ldots, |r_k - r_k'|_{p_k}^{a_k},
                                                |r_{k + 1} - r_{k + 1}'|)
\end{equation}
for each $r, r' \in R_S$.  It is easy to see that this is also a
translation-invariant metric on $R_S$ that determines the same topology
on $R_S$ as before.  This metric enjoys other properties like those in
the previous section as well.  More precisely, if $x, \in {\bf Z}_S$
and $x \ne y$, then
\begin{equation}
\label{d_{S, a}(Delta_S(x), Delta_S(y)) ge 1}
        d_{S, a}(\Delta_S(x), \Delta_S(y)) \ge 1,
\end{equation}
as in (\ref{d_S(Delta_S(x), Delta_S(y)) ge 1}).  Similarly, for each
$r \in R_S$, there is an $y \in {\bf Z}_S$ such that
\begin{equation}
\label{d_{S, a}(r, Delta_S(y)) le 1}
        d_{S, a}(r, \Delta_S(y)) \le 1,
\end{equation}
as in (\ref{d_S(r, Delta_S(y)) le 1}).  The proofs of these two
statements are essentially the same as before, because $|z|_{p_j} \le 1$
for some $z \in {\bf Q}_{p_j}$ implies that $|z|_{p_j}^{a_j} \le 1$ too.

        The Hausdorff dimension of $R_S$ with respect to
$d_{S, a}(r, r')$ is equal to
\begin{equation}
\label{1 + sum_{j = 1}^k 1/a_k}
        1 + \sum_{j = 1}^k 1/a_k,
\end{equation}
which reduces to $k + 1$ when $a_j = 1$ for each $j$.  This is the
expected Hausdorff dimension for a Cartesian product in nice
situations, and we shall say more about it in the next section.  The
topological dimension of $R_S$ is equal to $1$, because the
topological dimension of ${\bf R}$ is $1$, and the topological
dimension of ${\bf Q}_p$ is $0$ for each prime $p$.  Thus an advantage
of using the $a_j$'s is that the Hausdorff dimension can approximate
the topological dimension in a simple way when the $a_j$'s are
sufficiently large.

\section{Haar measure}
\label{haar measure}
\setcounter{equation}{0}

        Let us continue with the notation and assumptions in the
previous sections, especially Section \ref{finitely many primes}.
Let $H_S$ be the product measure on $R_S$ corresponding to Haar
measure $H_{p_j}$ on ${\bf Q}_{p_j}$ for $j = 1, \ldots, k$ and Lebesgue
measure on the real line, which is also Haar measure on ${\bf R}$.
This is a Haar measure on $R_S$ as a locally compact commutative
topological group with respect to addition, which is characterized by
the normalization
\begin{equation}
\label{H_S({bf Z}_{p_1} times cdots times {bf Z}_{p_k} times [0, 1]) = 1}
        H_S({\bf Z}_{p_1} \times \cdots \times {\bf Z}_{p_k} \times [0, 1]) = 1.
\end{equation}

        Suppose that $t \in R_S$ satisfies $t_j \ne 0$ for each
$j = 1, \ldots, k + 1$.  Thus the mapping $r \mapsto t \, r$ is a 
homeomorphism from $R_S$ onto itself, where $t \, r$ refers to
multiplication on $R_S$, which is defined coordinatewise.  If $E
\subseteq R_S$ is a Borel set, then it follows that
\begin{equation}
\label{t E = {t r : r in E}}
        t \, E = \{t \, r : r \in E\}
\end{equation}
is a Borel set too.  One can check that
\begin{equation}
\label{H_S(t E) = mu(t) H_S(E)}
        H_S(t \, E) = \mu(t) \, H_S(E),
\end{equation}
where
\begin{equation}
\label{mu(t) = (prod_{j = 1}^k |t_j|_{p_j}) |t_{k + 1}|}
        \mu(t) = \Big(\prod_{j = 1}^k |t_j|_{p_j}\Big) \, |t_{k + 1}|.
\end{equation}
This uses the corresponding properties of $H_{p_j}$ on ${\bf Q}_{p_j}$,
$1 \le j \le k$, and of Lebesgue measure on ${\bf R}$.

        In particular, one can apply this to $t = \Delta_S(x)$, where $x$
is a nonzero rational number.  Suppose that $x$ and $1/x$ are both in
${\bf Z}_S$, which happens exactly when
\begin{equation}
\label{x = pm p_1^{l_1} cdots p_k^{l_k}}
        x = \pm p_1^{l_1} \cdots p_k^{l_k}
\end{equation}
for some integers $l_1, \ldots, l_k$.  In this case, it is easy to see that
\begin{equation}
\label{mu(Delta_S(x)) = 1}
        \mu(\Delta_S(x)) = 1,
\end{equation}
so that
\begin{equation}
\label{H_S(Delta_S(x) E) = H_S(E)}
        H_S(\Delta_S(x) \, E) = H_S(E)
\end{equation}
for each Borel set $E \subseteq R_S$.

        Because the metric $d_S(r, r')$ on $R_S$ is defined in 
(\ref{d_S(r, r') = max(|r_1 - r'_1|_{p_1}, ..., |r_{k + 1} - r'_{k + 1}|)})
as the maximum of metrics in the $k + 1$ coordinates, a ball in $R_S$
of radius $\rho$ with respect to $d_S(r, r')$ is the same as the Cartesian 
product of balls of radius $\rho$ in the $k + 1$ factors.  It follows that
the $H_S$ measure of a ball of radius $\rho$ in $R_S$ with respect to
$d_S(r, r')$ is approximately $\rho^{k + 1}$, to within bounded factors,
because of the corresponding property of $H_{p_j}$ on ${\bf Q}_{p_j}$ for
$j = 1, \ldots, k$, and of Lebesgue measure on ${\bf R}$.  This implies
that $R_S$ is Ahlfors-regular of dimension $k + 1$ with respect to
$d_S(r, r')$, and in particular that $R_S$ has Hausdorff dimension $k + 1$
with respect to $d_S(r, r')$.  Similarly, if $d_{S, a}(r, r')$ is as in
(\ref{d_{S, a}(r, r') = max(...)}), then $H_S$ is Ahlfors-regular with
respect to $d_{S, a}(r, r')$, with dimension (\ref{1 + sum_{j = 1}^k 1/a_k}).
This implies that the Hausdorff dimension of $R_S$ with respect to
$d_{S, a}(r, r')$ is also equal to (\ref{1 + sum_{j = 1}^k 1/a_k}).

\section{Quotient spaces}
\label{quotient spaces}
\setcounter{equation}{0}

        As in Section \ref{finitely many primes}, $R_S$ is a locally
compact commutative topological ring with respect to coordinatewise
addition and multiplication, and $\Delta_S({\bf Z}_S)$ is a discrete
sub-ring in $R_S$.  In particular, $R_S$ is a locally compact
commutative topological group with respect to addition, and
$\Delta_S({\bf Z}_S)$ is a discrete subgroup of $R_S$.  Thus the quotient
$R_S / \Delta_S({\bf Z}_S)$ may be defined as a commutative topological
group in the usual way.  More precisely, $R_S / \Delta_S({\bf Z}_S)$ is
compact with respect to the quotient topology, because every element of
$R_S$ is at a bounded distance from an element of $\Delta_S({\bf Z}_S)$,
as before.

        Alternatively, consider the Cartesian product
\begin{equation}
\label{A_S = {bf Z}_{p_1} times cdots times {bf Z}_{p_k} times {bf R}}
        A_S = {\bf Z}_{p_1} \times \cdots \times {\bf Z}_{p_k} \times {\bf R},
\end{equation}
which is a sub-ring of $R_S$ that is both open and closed with respect
to the usual product topology on $R_S$.  Observe that $\Delta_S({\bf Z})
\subseteq A_S$, and in fact
\begin{equation}
\label{Delta_S({bf Z}) = Delta_S({bf Z}_S) cap A_S}
        \Delta_S({\bf Z}) = \Delta_S({\bf Z}_S) \cap A_S,
\end{equation}
which is the same as saying that $x \in {\bf Z}_S$ is an integer if
and only if $|x|_{p_j} \le 1$ for each $j = 1, \ldots, k$.  As before,
$\Delta_S({\bf Z})$ is a discrete subgroup of $A_S$ as a topological
group with respect to addition, so that the quotient $A_S /
\Delta_S({\bf Z})$ may be defined as a commutative topological group.
It is even easier to see that every element of $A_S$ is at a bounded
distance from an element of $\Delta_S({\bf Z})$ in this case, and
hence that $A_S / \Delta_S({\bf Z})$ is compact.

        The obvious inclusion of $A_S$ in $R_S$ leads to a natural
group homomorphism from $A_S$ into $R_S / \Delta_S({\bf Z}_S)$, and
the kernel of this homomorphism is equal to $\Delta_S({\bf Z})$, by
(\ref{Delta_S({bf Z}) = Delta_S({bf Z}_S) cap A_S}).  Thus we get a
natural group homomorphism from $A_S / \Delta_S({\bf Z})$ into $R_S /
\Delta_S({\bf Z}_S)$, whose kernel is trivial.  One can check that
these two natural homomorphisms are surjective, by an argument like
the one used in Section \ref{finitely many primes} to show that every
element of $R_S$ is at bounded distance from an element of
$\Delta_S({\bf Z}_S)$.  More precisely, for each $r \in R_S$, one can
choose $x \in {\bf Z}_S$ as in (\ref{x = sum_{j = 1}^k x_j}), to get that
\begin{equation}
\label{r - Delta_S(x) in A_S}
        r - \Delta_S(x) \in A_S,
\end{equation}
as desired.  It is easy to see that this natural isomorphism between
$A_S / \Delta_S({\bf Z})$ and $R_S / \Delta_S({\bf Z}_S)$ is also a
homeomorphism with respect to the corresponding quotient topologies,
which are locally the same as the product topology on $R_S$.

        Suppose that $x \in {\bf Z}_S$, $x \ne 0$, and $1/x \in {\bf Z}_S$,
which is equivalent to the representation (\ref{x = pm p_1^{l_1} cdots
  p_k^{l_k}}) in the previous section.  As before,
\begin{equation}
\label{r mapsto Delta_S(x) r}
        r \mapsto \Delta_S(x) \, r
\end{equation}
is a homeomorphism from $R_S$ onto itself, where the product on the
right uses coordinatewise multiplication on $R_S$.  Of course, (\ref{r
  mapsto Delta_S(x) r}) is also an automorphism on $R_S$ as a
commutative group with respect to addition.  Note that (\ref{r mapsto
  Delta_S(x) r}) maps $\Delta_S({\bf Z}_S)$ onto itself under these
conditions as well.  It follows that (\ref{r mapsto Delta_S(x) r})
induces an automorphism on $R_S / \Delta_S({\bf Z}_S)$ as a
commutative topological group with respect to the usual quotient
topology.

\section{Connectedness}
\label{connectedness}
\setcounter{equation}{0}

        Consider the Cartesian product
\begin{equation}
\label{{bf Q}_{p_1} times cdots times {bf Q}_{p_k}}
        {\bf Q}_{p_1} \times \cdots \times {\bf Q}_{p_k},
\end{equation}
equipped with the product topology associated to the usual topology on
${\bf Q}_{p_j}$ for each $j = 1, \ldots, k$.  This is a locally
compact commutative topological ring with respect to coordinatewise
addition and multiplication.  Also let $\delta_S$ be the diagonal
embedding of ${\bf Q}$ into (\ref{{bf Q}_{p_1} times cdots times {bf
    Q}_{p_k}}), which sends each $x \in {\bf Q}$ to the $k$-tuple
whose coordinates are all equal to $x$.  Thus $\delta_S$ is a ring
homomorphism from ${\bf Q}$ into (\ref{{bf Q}_{p_1} times cdots times
  {bf Q}_{p_k}}).  Note that
\begin{equation}
\label{{bf Z}_{p_1} times cdots times {bf Z}_{p_k}}
        {\bf Z}_{p_1} \times \cdots \times {\bf Z}_{p_k}
\end{equation}
is a compact open sub-ring of (\ref{{bf Q}_{p_1} times cdots times {bf
    Q}_{p_k}}).

        Of course, $\delta_S({\bf Z})$ is contained in 
(\ref{{bf Z}_{p_1} times cdots times {bf Z}_{p_k}}), and in fact the 
closure of $\delta_S({\bf Z})$ in (\ref{{bf Q}_{p_1} times cdots times
 {bf Q}_{p_k}}) is equal to (\ref{{bf Z}_{p_1} times cdots times {bf Z}_{p_k}}).
To see this, it suffices to check that the closure of $\delta_S({\bf
  Z})$ in (\ref{{bf Q}_{p_1} times cdots times {bf Q}_{p_k}}) contains
${\bf Z}^k$, because ${\bf Z}$ is dense in ${\bf Z}_p$ for each prime
$p$.  If $n$ is any positive integer, and if $z_1, \ldots, z_k$ are
integers, then it is well known that there is an integer $z$ such that
$z \equiv z_k$ modulo $p_k^n$, because $p_1, \ldots, p_k$ are distinct
primes.  This says exactly that $(z_1, \ldots, z_k)$ can be
approximated by $\delta_S(z)$ in (\ref{{bf Q}_{p_1} times cdots times
  {bf Q}_{p_k}}), as desired.  It follows that $\delta_S({\bf Z}_S)$
is dense in (\ref{{bf Q}_{p_1} times cdots times {bf Q}_{p_k}}), and
in particular that $\delta_S({\bf Q})$ is dense in (\ref{{bf Q}_{p_1}
  times cdots times {bf Q}_{p_k}}).

        There is an obvious isomorphism between the real line and
\begin{equation}
\label{{0} times cdots times {0} times {bf R}}
        \{0\} \times \cdots \times \{0\} \times {\bf R},
\end{equation}
as a closed subring of $R_S$.  It is easy to see that the natural
quotient mapping from $R_S$ onto $R_S / \Delta_S({\bf Z}_S)$ maps
(\ref{{0} times cdots times {0} times {bf R}}) to a dense subset of
$R_S / \Delta_S({\bf Z}_S)$, since $\delta_S({\bf Z}_S)$ is dense
in (\ref{{bf Q}_{p_1} times cdots times {bf Q}_{p_k}}).  Equivalently,
the natural quotient mapping from $A_S$ onto $A_S / \Delta({\bf Z})$
maps (\ref{{0} times cdots times {0} times {bf R}}) to a dense subset of
$A_S / \Delta_S({\bf Z})$, because $\delta_S({\bf Z})$ is dense in
(\ref{{bf Z}_{p_1} times cdots times {bf Z}_{p_k}}).  It follows that
$R_S / \Delta_S({\bf Z}_S) \cong A_S / \Delta_S({\bf Z})$ is connected,
because the real line is connected, and the closure of a connected set
is connected.

        There is a natural continuous ring homomorphism from $R_S$
onto ${\bf R}$, which sends $r \in R_S$ to $r_{k + 1} \in {\bf R}$.
Of course, ${\bf Z}$ is a discrete subgroup of ${\bf R}$ as a locally
compact commutative topological group with respect to addition, and
the quotient ${\bf R} / {\bf Z}$ is isomorphic as a topological group
to the multiplicative group of complex numbers with modulus equal to $1$.
Consider the mapping from $A_S$ onto ${\bf R} / {\bf Z}$ which sends
$r \in A_S$ to the image of $r_{k + 1}$ in ${\bf R} / {\bf Z}$.  This
is a continuous group homomorphism from $A_S$ as a topological group
with respect to addition onto ${\bf R} / {\bf Z}$, and the kernel of
this homomorphism is
\begin{equation}
\label{{bf Z}_{p_1} times cdots times {bf Z}_{p_k} times {bf Z}}
        {\bf Z}_{p_1} \times \cdots \times {\bf Z}_{p_k} \times {\bf Z}.
\end{equation}
In particular, (\ref{{bf Z}_{p_1} times cdots times {bf Z}_{p_k} times
  {bf Z}}) contains $\Delta_S({\bf Z})$ as a subgroup, and so this
homomorphism from $A_S$ onto ${\bf R} / {\bf Z}$ leads to a continuous
group homomorphism from $A_S / \Delta_S({\bf Z})$ onto ${\bf R} / {\bf
  Z}$.  Note that the kernel of this homomorphism from $A_S /
\Delta_S({\bf Z})$ onto ${\bf R} / {\bf Z}$ corresponds to the
quotient of (\ref{{bf Z}_{p_1} times cdots times {bf Z}_{p_k} times
  {bf Z}}) by $\Delta_S({\bf Z})$.  Using this, it is easy to see that
the kernel of this homomorphism from $A_S / \Delta_S({\bf Z})$ onto
${\bf R} / {\bf Z}$ is isomorphic to (\ref{{bf Z}_{p_1} times cdots
  times {bf Z}_{p_k}}) as a commutative topological group with respect
to addition.

\section{Quotient metrics}
\label{quotient metrics}
\setcounter{equation}{0}

        Let $q_S$ be the natural quotient mapping from $R_S$ onto
$R_S / \Delta_S({\bf Z}_S)$.  Also let $a = (a_1, \ldots, a_k)$ be a
$k$-tuple of positive real numbers, and let $d_{S, a}(r, t)$ be the
corresponding metric on $R_S$ defined in (\ref{d_{S, a}(r, r') = max(...)}).
This leads to a quotient metric on $R_S / \Delta_S({\bf Z}_S)$, given by
\begin{eqnarray}
\label{D_{S, a}(q_s(r), q_s(t)) = ...}
\lefteqn{D_{S, a}(q_S(r), q_S(t))}                               \\
 & = & \inf \{d_{S, a}(\widetilde{r}, \widetilde{t}) :
                    \widetilde{r}, \widetilde{t} \in R_S,
   \, q_S(\widetilde{r}) = q_S(r), \, q_S(\widetilde{t}) = q_S(t)\}. \nonumber
\end{eqnarray}
Equivalently,
\begin{eqnarray}
\label{D_{S, a}(q_S(r), q_S(t)) = ..., 2}
 D_{S, a}(q_S(r), q_S(t)) & = & \inf \{d_{S, a}(\widetilde{r}, t) :
                     \widetilde{r} \in R_S, \, q_S(\widetilde{r}) = q_S(r)\} \\
                         & = & \inf \{d_{S, a}(r, \widetilde{t}) : 
              \widetilde{t} \in R_S, \, q_S(\widetilde{t}) = q_S(t)\} \nonumber
\end{eqnarray}
for every $r, t \in R_S$, because $d_{S, a}(\cdot, \cdot)$ is
invariant under translations on $R_S$.  In order to check that $D_{S,
  a}(\cdot, \cdot)$ satisfies the triangle inequality on $R_S /
\Delta_S({\bf Z}_S)$, let $r, t, u \in R_S$ be given, and suppose that
$\widetilde{r}, \widetilde{u} \in R_S$ satisfy $q_S(\widetilde{r}) =
q_S(r)$, $q_S(\widetilde{u}) = q_S(u)$.  Using the definition
(\ref{D_{S, a}(q_s(r), q_s(t)) = ...}) of $D_{S, a}(q_S(r), q_S(u))$
and then the triangle inequality for $d_{S, a}(\cdot, \cdot)$, we get that
\begin{equation}
\label{D_{S, a}(q_S(r), q_S(u)) le ...}
        D_{S, a}(q_S(r), q_S(u)) \le d_{S, a}(\widetilde{r}, \widetilde{u}) 
                 \le d_{S, a}(\widetilde{r}, t) + d_{S, a}(t, \widetilde{u}).
\end{equation}
Taking the infimum over $\widetilde{r}$ and $\widetilde{u}$ and using
(\ref{D_{S, a}(q_S(r), q_S(t)) = ..., 2}), we get that
\begin{equation}
\label{D_{S, a}(q_S(r), q_S(u)) le ..., 2}
 D_{S, a}(q_S(r), q_S(u)) \le D_{S, a}(q_S(r), q_S(t)) + D_{S, a}(q_S(t), q_S(u)),
\end{equation}
as desired.  It is easy to see that $D_{S, a}(q_S(r), q_S(t))$ is
symmetric in $r$, $t$ and invariant under translations, because of the
analogous properties of $d_{S, a}(r, t)$.  One can also verify that
$D_{S, a}(q_S(r), q_S(t)) = 0$ if and only if $q_S(r) = q_S(t)$ in
$R_S / \Delta_S({\bf Z}_S)$, basically because $\Delta_S({\bf Z}_S)$
is a closed subgroup of $R_S$.  This shows that $D_{S, a}(\cdot, \cdot)$
is indeed a metric on $R_S / \Delta_S({\bf Z}_S)$.

        Similarly, let $q_S'$ be the natural quotient mapping from $A_S$
onto $A_S / \Delta_S({\bf Z})$, and define $D_{S, a}'(\cdot, \cdot)$
on $A_S / \Delta_S({\bf Z})$ by
\begin{eqnarray}
\label{D_{S, a}'(q_S'(r), q_S'(t)) = ...}
\lefteqn{D_{S, a}'(q_S'(r), q_S'(t))}                       \\
 & = & \inf \{d_{S, a}(\widetilde{r}, \widetilde{t}) :
                                     \widetilde{r}, \widetilde{t} \in A_S,
\, q_S'(\widetilde{r}) = q_S'(r), \, q_S'(\widetilde{t}) = q_S'(t)\}. \nonumber
\end{eqnarray}
As before, this is the same as
\begin{eqnarray}
\label{D_{S, a}'(q_S'(r), q'_S(t)) = ..., 2}
 D_{S, a}'(q_S'(r), q'_S(t)) & = & \inf \{d_{S, a}(\widetilde{r}, t) :
                   \widetilde{r} \in A_S, \, q_S'(\widetilde{r}) = q_S'(r)\} \\
                            & = & \inf \{d_{S, a}(r, \widetilde{t}) :
            \widetilde{t} \in A_S, \, q_S'(\widetilde{t}) = q_S'(t)\} \nonumber
\end{eqnarray}
for every $r, t \in A_S$, because $d_{S, a}(\cdot, \cdot)$ is
invariant under translations.  One can use this to check that $D_{S,
  a}'(\cdot, \cdot)$ satisfies the triangle inequality
\begin{equation}
\label{D_{S, a}'(q_S'(r), q_S'(u)) le ...}
        D_{S, a}'(q_S'(r), q_S'(u)) \le D_{S, a}'(q_S'(r), q_S'(t)) 
                                           + D_{S, a}'(q_S'(t), q_S'(u))
\end{equation}
for every $r, t, u \in A_S$, as in the previous situation.  Clearly
$D_{S, a}'(q_S'(r), q_S'(t))$ is symmetric in $r$, $t$ and invariant
under translations, because of the corresponding properties of $d_{S,
  a}(\cdot, \cdot)$.  One can use the fact that $\Delta_S({\bf Z})$ is
a closed subgroup of $A_S$ to show that $D_{S, a}'(q_S'(r), q_S'(t)) =
0$ if and only if $q_S'(r) = q_S'(t)$ in $A_S / \Delta_S({\bf Z})$,
which implies that $D_{S, a}'(\cdot, \cdot)$ is a metric on $A_S /
\Delta_S({\bf Z})$.

        It is easy to see that
\begin{equation}
\label{D_{S, a}'(q_S'(r), q_S'(t)) le 1}
        D_{S, a}'(q_S'(r), q_S'(t)) \le 1
\end{equation}
for every $r, t \in A_S$.  Of course, $|x - y|_p \le 1$ for every $x,
y \in {\bf Z}_p$ and for each prime $p$, and so one only has to
consider the contribution of the standard metric on ${\bf R}$ to
$d_{S, a}(\cdot, \cdot)$ in the definition of $D_{S, a}'(q_S'(r),
q_S'(t))$ to get (\ref{D_{S, a}'(q_S'(r), q_S'(t)) le 1}).  Thus
(\ref{D_{S, a}'(q_S'(r), q_S'(t)) le 1}) reduces to the fact that
there is a $z \in {\bf Z}$ such that
\begin{equation}
\label{|r_{k + 1} - t_{k + 1} - z| le 1/2}
        |r_{k + 1} - t_{k + 1} - z| \le 1/2.
\end{equation}

        Next, let us check that
\begin{equation}
\label{D_{S, a}(q_S(r), q_S(t)) le D_{S, a}'(q_S'(r), q_S'(t))}
        D_{S, a}(q_S(r), q_S(t)) \le D_{S, a}'(q_S'(r), q_S'(t))
\end{equation}
for every $r, t \in A_S$.  More precisely, if $r, \widetilde{r}, t,
\widetilde{t} \in A_S$ satisfy $q_S'(\widetilde{r}) = q_S'(r)$ and
$q_S'(\widetilde{t}) = q_S'(t)$ in $A_S / \Delta_S({\bf Z})$, then
$q_S(\widetilde{r}) = q_S(r)$ and $q_S(\widetilde{t}) = q_S(t)$
in $R_S / \Delta_S({\bf Z}_S)$, because ${\bf Z} \subseteq {\bf Z}_S$.
Hence the distances $d_{S, a}(\widetilde{r}, \widetilde{t})$ in the
infimum in (\ref{D_{S, a}'(q_S'(r), q_S'(t)) = ...}) are also included
in the infimum in (\ref{D_{S, a}(q_s(r), q_s(t)) = ...}), which implies
(\ref{D_{S, a}(q_S(r), q_S(t)) le D_{S, a}'(q_S'(r), q_S'(t))}).

        Suppose now that $r, t \in A_S$, $\widetilde{t} \in R_S$, and
$q_S(\widetilde{t}) = q_S(t)$, and let us check that
\begin{equation}
\label{D_{S, a}'(q_S'(r), q_S'(t)) le d_{S, a}(r, widetilde{t})}
        D_{S, a}'(q_S'(r), q_S'(t)) \le d_{S, a}(r, \widetilde{t}).
\end{equation}
If $\widetilde{t} \in A_S$, then $\widetilde{t} - t \in A_S \cap
\Delta_S({\bf Z}_S) = \Delta_S({\bf Z})$, as in (\ref{Delta_S({bf Z})
  = Delta_S({bf Z}_S) cap A_S}).  This implies that
$q_S'(\widetilde{t}) = q_S'(t)$, in which case (\ref{D_{S,
    a}'(q_S'(r), q_S'(t)) le d_{S, a}(r, widetilde{t})}) follows from
(\ref{D_{S, a}'(q_S'(r), q'_S(t)) = ..., 2}).  Otherwise, if
$\widetilde{t} \not\in A_S$, then $\widetilde{t}_j \not\in {\bf Z}_{p_j}$
for some $j = 1, \ldots, k$, which means that $|\widetilde{t}_j|_{p_j} > 1$.
However, $r_j \in {\bf Z}_{p_j}$, because $r \in A_S$, and
\begin{equation}
\label{|widetilde{t}_j|_{p_j} le max(|r_j - widetilde{t}_j|_{p_j}, |r_j|_{p_j})}
 |\widetilde{t}_j|_{p_j} \le \max(|r_j - \widetilde{t}_j|_{p_j}, |r_j|_{p_j}),
\end{equation}
by the ultrametric version (\ref{|x + y|_p le max(|x|_p, |y|_p)}) of
the triangle inequality.  Thus we get
\begin{equation}
\label{|r_j - widetilde{t}_j|_{p_j} > 1}
        |r_j - \widetilde{t}_j|_{p_j} > 1,
\end{equation}
since $|\widetilde{t}_j|_{p_j} > 1$ and $|r_j|_{p_j} \le 1$.  This
shows that
\begin{equation}
\label{d_{S, a}(r, widetilde{t}) ge |r_j - widetilde{t}_j|_{p_j}^{a_j} > 1}
        d_{S, a}(r, \widetilde{t}) \ge |r_j - \widetilde{t}_j|_{p_j}^{a_j} > 1
\end{equation}
in this case, which implies (\ref{D_{S, a}'(q_S'(r), q_S'(t)) le d_{S,
    a}(r, widetilde{t})}), because of (\ref{D_{S, a}'(q_S'(r),
  q_S'(t)) le 1}).

        It follows that
\begin{equation}
\label{D_{S, a}'(q_S'(r), q_S'(t)) le D_{S, a}(q_S(r), q_S(t))}
        D_{S, a}'(q_S'(r), q_S'(t)) \le D_{S, a}(q_S(r), q_S(t))
\end{equation}
for every $r, t \in A_S$, by (\ref{D_{S, a}(q_S(r), q_S(t)) = ..., 2})
and (\ref{D_{S, a}'(q_S'(r), q_S'(t)) le d_{S, a}(r, widetilde{t})}).
Combining this with (\ref{D_{S, a}(q_S(r), q_S(t)) le D_{S,
    a}'(q_S'(r), q_S'(t))}), we get that
\begin{equation}
        D_{S, a}'(q_S'(r), q_S'(t)) = D_{S, a}(q_S(r), q_S(t))
\end{equation}
for every $r, t \in A_S$.  Thus the metrics $D_{S, a}(\cdot, \cdot)$
on $R_S / \Delta_S({\bf Z}_S)$ and $D_{S, a}'(\cdot, \cdot)$ on $A_S /
\Delta_S({\bf Z})$ correspond exactly to each other under the natural
isomorphism between $R_S / \Delta_S({\bf Z}_S)$ and $A_S /
\Delta_S({\bf Z})$ described in Section \ref{quotient spaces}.

        Of course,
\begin{equation}
\label{D_{S, a}'(q_S'(r), q_S'(t)) le d_{S, a}(r, t)}
        D_{S, a}'(q_S'(r), q_S'(t)) \le d_{S, a}(r, t)
\end{equation}
for every $r, t \in A_S$, by the definition (\ref{D_{S, a}'(q_S'(r),
  q_S'(t)) = ...}) of $D_{S, a}'(q_S'(r), q_S'(t))$.  Suppose that
\begin{equation}
\label{|r_{k + 1} - t_{k + 1}| le 1/2}
        |r_{k + 1} - t_{k + 1}| \le 1/2,
\end{equation}
which implies that
\begin{equation}
\label{|r_{k + 1} - t_{k + 1} - z| ge 1/2}
        |r_{k + 1} - t_{k + 1} - z| \ge 1/2
\end{equation}
for every nonzero integer $z$.  If $\widetilde{t} \in A_S$,
$q_S'(\widetilde{t}) = q_S'(t)$, and $\widetilde{t} \ne t$, then
\begin{equation}
\label{widetilde{t}_{k + 1} = t_{k + 1} + z}
        \widetilde{t}_{k + 1} = t_{k + 1} + z
\end{equation}
for some nonzero integer $z$, and hence
\begin{equation}
\label{d_{S, a}(r, widetilde{t}) ge |r_{k + 1} - widetilde{t}_{k + 1}| ge 1/2}
 d_{S, a}(r, \widetilde{t}) \ge |r_{k + 1} - \widetilde{t}_{k + 1}| \ge 1/2.
\end{equation}
This implies that
\begin{equation}
\label{D_{S, a}'(q_S'(r), q_S'(t)) ge min(d_{S, a}(r, t), 1/2)}
        D_{S, a}'(q_S'(r), q_S'(t)) \ge \min(d_{S, a}(r, t), 1/2),
\end{equation}
by considering the cases where $\widetilde{t} = t$ and $\widetilde{t}
\ne t$ in (\ref{D_{S, a}'(q_S'(r), q'_S(t)) = ..., 2}).  Clearly
\begin{equation}
\label{d_{S, a}(r, t) le 1}
        d_{S, a}(r, t) \le 1
\end{equation}
when $r, t \in A_S$ satisfy (\ref{|r_{k + 1} - t_{k + 1}| le 1/2}),
and so (\ref{D_{S, a}'(q_S'(r), q_S'(t)) ge min(d_{S, a}(r, t), 1/2)})
implies that
\begin{equation}
\label{d_{S, a}(r, t) le 2 D_{S, a}'(q_S'(r), q_S'(t))}
        d_{S, a}(r, t) \le 2 \, D_{S, a}'(q_S'(r), q_S'(t))
\end{equation}
under these conditions.  Alternatively, if $d_{S, a}(r, t) \le 1/2$,
then (\ref{|r_{k + 1} - t_{k + 1}| le 1/2}) holds in particular, and
we can combine (\ref{D_{S, a}'(q_S'(r), q_S'(t)) le d_{S, a}(r, t)})
and (\ref{D_{S, a}'(q_S'(r), q_S'(t)) ge min(d_{S, a}(r, t), 1/2)}) to
get that
\begin{equation}
\label{D_{S, a}'(q_S'(r), q_S'(t)) = d_{S, a}(r, t)}
        D_{S, a}'(q_S'(r), q_S'(t)) = d_{S, a}(r, t)
\end{equation}
in this case.

\section{$n$ Dimensions}
\label{n dimensions}
\setcounter{equation}{0}

        Let $n$ be a positive integer, and let $R$ be a commutative ring
with nonzero multiplicative identity element $e$.  As usual, the set
$M_n(R)$ of $n \times n$ matrices with entries in $R$ is a ring with
respect to matrix addition and multiplication.  The matrix $I$ with
diagonal entries equal to $e$ and off-diagonal entries equal to $0$ is
the identity element in $M_n(R)$, and the group of invertible matrices
in $M_n(R)$ is denoted $GL(n, R)$.  By standard arguments, an element
of $M_n(R)$ is invertible if and only if its determinant is an
invertible element of $R$.  The set of matrices with determinant equal
to $e$ is denoted $SL(n, R)$, which is a subgroup of $GL(n, R)$.

        Because ${\bf Z}_S$ is a sub-ring of the field ${\bf R}$ of real 
numbers, there is a natural embedding of $SL(n, {\bf Z}_S)$ into
$SL(n, {\bf R})$.  Similarly, ${\bf Z}_S$ may be considered as a sub-ring 
of ${\bf Q}_{p_j}$ for $j = 1, \ldots, k$, which leads to a natural
embedding of $SL(n, {\bf Z}_S)$ into $SL(n, {\bf Q}_{p_j})$.  The
combination of these embeddings leads to a diagonal embedding of
$SL(n, {\bf Z}_S)$ onto a discrete subgroup of
\begin{equation}
\label{SL(n, Q_{p_1}) times cdots times SL(n, Q_{p_k}) times SL(n, R)}
        SL(n, {\bf Q}_{p_1}) \times \cdots \times SL(n, {\bf Q}_{p_k})
                                                         \times SL(n, {\bf R}).
\end{equation}
Alternatively, one can use the embedding $\Delta_S$ of ${\bf Z}_S$
into $R_S$ to embed $SL(n, {\bf Z}_S)$ onto a discrete subgroup of
$SL(n, R_S)$, and to embed $GL(n, {\bf Z}_S)$ onto a discrete subgroup
of $GL(n, R_S)$.  More precisely, one can use $\Delta_S$ to embed
$M_n({\bf Z}_S)$ onto a discrete sub-ring of $M_n(R_S)$, and this
embedding intertwines determinants on $M_n({\bf Z}_S)$ and $M_n(R_S)$
in the obvious way.

        Of course, $GL(n, {\bf Z})$ is already a discrete subgroup of
$GL(n, {\bf R})$.  Elements of $GL(n, {\bf R})$ determine invertible 
linear mappings on ${\bf R}^n$ in the usual way, and the linear
mappings on ${\bf R}^n$ corresponding to elements of $GL(n, {\bf Z})$
send ${\bf Z}^n$ onto itself.  Hence these mappings induce invertible
mappings on the $n$-dimensional torus ${\bf R}^n / {\bf Z}^n$, which
defines a well-known action of $GL(n, {\bf Z})$  on ${\bf R}^n / {\bf Z}^n$.

        Similarly, elements of $GL(n, R_S)$ determine invertible mappings 
on $R_S^n$, and we can use $\Delta_S$ to embed $GL(n, {\bf Z}_S)$ into
$GL(n, R_S)$, as before.  The mappings on $R_S^n$ associated to
elements of $GL(n, {\bf Z}_S)$ in this way send $\Delta_S({\bf
  Z}_S)^n$ onto itself, and thus induce invertible mappings on
\begin{equation}
\label{R_S^n / Delta_S({bf Z}_S)^n cong (R_S / Delta_S({bf Z}_S))^n}
        R_S^n / \Delta_S({\bf Z}_S)^n \cong (R_S / \Delta_S({\bf Z}_S))^n.
\end{equation}
 This leads to a natural action of $GL(n, {\bf Z}_S)$ on (\ref{R_S^n /
   Delta_S({bf Z}_S)^n cong (R_S / Delta_S({bf Z}_S))^n}), which is
 analogous to the classical action of $GL(n, {\bf Z})$ on ${\bf R}^n /
 {\bf T}^n$.

        If $T$ is any element of $M_n({\bf R})$, then $T$ still determines 
a linear mapping $\widehat{T}$ on ${\bf R}^n$ in the usual way.  If $T
\in M_n({\bf Z})$, then $\widehat{T}({\bf Z}^n) \subseteq {\bf Z}^n$,
and hence $\widehat{T}$ induces a mapping from ${\bf R}^n / {\bf Z}^n$
into itself.  If
\begin{equation}
\label{T in GL(n, {bf R}) cap M_n({bf Z})}
        T \in GL(n, {\bf R}) \cap M_n({\bf Z}),
\end{equation}
then the induced mapping on ${\bf R}^n / {\bf Z}^n$ is surjective, but
may not be one-to-one.  In this case, the determinant $\det T$ of $T$
is a nonzero integer.  If each prime factor of $\det T$ is one of the
$p_j$'s, then $\det T$ is an invertible element of ${\bf Z}_S$, so
that $T \in GL(n, {\bf Z}_S)$, and the remarks in the previous
paragraph apply.

        As before, we can use $\Delta_S$ to embed $M_n({\bf Z}_S)$ onto a
discrete sub-ring of $M_n(R_S)$.  Elements of $M_n(R_S)$ determine
mappings on $R_S^n$ in the usual way, and the mappings on $R_S^n$
associated to elements of $M_n({\bf Z}_S)$ map $\Delta_S({\bf Z}_S)^n$
into itself.  Thus each $T \in M_n({\bf Z}_S)$ leads to a mapping from
(\ref{R_S^n / Delta_S({bf Z}_S)^n cong (R_S / Delta_S({bf Z}_S))^n})
into itself.

        Note that $\Delta_S(x)$ is invertible in $R_S$ for every
$x \in {\bf Q}$ with $x \ne 0$.  If $T$ is an element of $M_n({\bf Z}_S)$ 
such that $\det T \ne 0$, then it follows that $T$ corresponds to an
element of $GL(n, R_S)$ under $\Delta_S$.  This implies that the
associated mapping on $R_S^n$ is invertible, and that the induced
mapping on (\ref{R_S^n / Delta_S({bf Z}_S)^n cong (R_S / Delta_S({bf
    Z}_S))^n}) is surjective.  However, $\det T$ may not be invertible
in ${\bf Z}_S$, even though it is a nonzero element of ${\bf Z}_S$.
Hence the associated mapping on $R_S^n$ may not send $\Delta_S({\bf
  Z}_S)^n$ onto itself, and the induced mapping on (\ref{R_S^n /
  Delta_S({bf Z}_S)^n cong (R_S / Delta_S({bf Z}_S))^n}) may not be
one-to-one.  This is analogous to (\ref{T in GL(n, {bf R}) cap M_n({bf
    Z})}) for mappings on ${\bf R}^n / {\bf Z}^n$.

\section{Some approximations}
\label{some approximations}
\setcounter{equation}{0}

        Let $l = (l_1, \ldots, l_k)$ be a $k$-tuple of nonnegative integers,
and put
\begin{equation}
\label{p^l = p_1^{l_1} cdots p_k^{l_k}}
        p^l = p_1^{l_1} \cdots p_k^{l_k},
\end{equation}
which is a positive integer.  Also put
\begin{equation}
\label{B_l = (Z / p_1^{l_1} Z) times cdots times (Z / p_k^{l_k} Z) times R}
        B_l = ({\bf Z} / p_1^{l_1} \, {\bf Z}) \times \cdots \times
                         ({\bf Z} / p_k^{l_k} \, {\bf Z}) \times {\bf R},
\end{equation}
which is a commutative topological ring with respect to coordinatewise
addition and multiplication, and using the product topology associated
to the standard topology on ${\bf R}$ and the discrete topology on
${\bf Z} / p_j^{l_j} \, {\bf Z}$ for $j = 1, \ldots, k$.  Remember
that there is a natural ring isomorphism between ${\bf Z}_{p_j} /
p_j^{l_j} \, {\bf Z}_{p_j}$ and ${\bf Z} / p_j^{l_j} \, {\bf Z}$ for each $j$,
which leads to a natural surjective ring homomorphism
\begin{equation}
\label{pi_{l, j} : {bf Z}_{p_j} to {bf Z} / p_j^{l_j} {bf Z}}
        \pi_{l, j} : {\bf Z}_{p_j} \to {\bf Z} / p_j^{l_j} \, {\bf Z}
\end{equation}
with kernel $p_j^{l_j} \, {\bf Z}_{p_j}$.  Thus
\begin{equation}
\label{pi_l(r) = (pi_{l, 1}(r_1), ldots, pi_{l, k}(r_k), r_{k + 1})}
        \pi_l(r) = (\pi_{l, 1}(r_1), \ldots, \pi_{l, k}(r_k), r_{k + 1})
\end{equation}
defines a continuous surjective ring homomorphism from $A_S$ onto $B_l$,
with kernel
\begin{equation}
\label{(p_1^{l_1} Z_{p_1}) times cdots times (p_k^{l_k} Z_{p_k}) times {0}}
        (p_1^{l_1} \, {\bf Z}_{p_1}) \times \cdots \times 
                             (p_k^{l_k} \, {\bf Z}_{p_k}) \times \{0\}.
\end{equation}

        The composition of $\pi_l$ with the diagonal embedding $\Delta_S$
of ${\bf Z}$ into $A_S$ defines an embedding of ${\bf Z}$ onto a discrete
sub-ring of $B_l$, and we let
\begin{equation}
\label{C_l = B_l / pi_l(Delta_S({bf Z}))}
        C_l = B_l / \pi_l(\Delta_S({\bf Z}))
\end{equation}
be the corresponding quotient group with respect to addition.  More
precisely, $C_l$ is a compact commutative topological group with
respect to the usual quotient topology associated to the topology
already defined on $B_l$.  The composition of $\pi_l$ with the standard
quotient mapping from $B_l$ onto $C_l$ defines a group homomorphism
from $A_S$ onto $C_l$ with respect to addition, and the kernel of
this homomorphism contains $\Delta_S({\bf Z})$ by construction.
This leads to a natural continuous surjective group homomorphism
\begin{equation}
\label{pi_l' : A_S / Delta_S({bf Z}) to C_l}
        \pi_l' : A_S / \Delta_S({\bf Z}) \to C_l,
\end{equation}
whose kernel is the image of (\ref{(p_1^{l_1} Z_{p_1}) times cdots
  times (p_k^{l_k} Z_{p_k}) times {0}}) under the natural quotient
mapping $q_S'$ from $A_S$ onto $A_S / \Delta_S({\bf Z})$.  Note that
the restriction of $q_S'$ to (\ref{(p_1^{l_1} Z_{p_1}) times cdots
  times (p_k^{l_k} Z_{p_k}) times {0}}) is injective, since the
intersection of (\ref{(p_1^{l_1} Z_{p_1}) times cdots times (p_k^{l_k}
  Z_{p_k}) times {0}}) with $\Delta_S({\bf Z})$ is trivial.

        There is a very simple way to embed ${\bf Z} \times {\bf R}$
into $A_S$, using the natural embedding of ${\bf Z}$ in ${\bf Z}_{p_j}$
for each $j = 1, \ldots, k$.  Equivalently, one combines the
diagonal embedding $\delta_S$ of ${\bf Z}$ in ${\bf Z}_{p_1} \times 
\cdots \times {\bf Z}_{p_k}$ with the identity mapping on ${\bf R}$,
to send $(x, y) \in {\bf Z} \times {\bf R}$ to $(x, \ldots, x, y) \in A_S$.
The image of ${\bf Z} \times {\bf R}$ in $A_S$ is dense in $A_S$,
because $\delta_S({\bf Z})$ is dense in ${\bf Z}_{p_1} \times \cdots \times
{\bf Z}_{p_k}$, as mentioned previously. 

        The composition of this embedding with the natural homomorphism 
$\pi_l$ from $A_S$ onto $B_l$ is also easy to describe.  It sends 
$(x, y) \in {\bf Z} \times {\bf R}$ to the point in  $B_l$ whose $j$th 
coordinate is the image of $x$ under the usual quotient mapping from 
${\bf Z}$ onto ${\bf Z} / p_j^{l_j} \, {\bf Z}$ for $j = 1, \ldots, k$, 
and whose $(k + 1)$rst coordinate is $y$.  This is a ring homomorphism
from ${\bf Z} \times {\bf R}$ onto $B_l$, whose kernel is $(p^l \, {\bf Z})
\times \{0\}$.

        Consider the diagonal mapping $\zeta : {\bf Z} \to {\bf Z} \times
{\bf R}$, which is defined by $\zeta(x) = (x, x)$.  The composition of 
$\zeta$ with the obvious embedding of ${\bf Z} \times {\bf R}$ into $A_S$
is the same as the diagonal mapping $\Delta_S$ from ${\bf Z}$ into $A_S$.
This leads to a natural embedding of
\begin{equation}
\label{({bf Z} times {bf R}) / zeta({bf Z})}
         ({\bf Z} \times {\bf R}) / \zeta({\bf Z})
\end{equation}
onto a dense subgroup of $A_S / \Delta_S({\bf Z})$.  Note that
(\ref{({bf Z} times {bf R}) / zeta({bf Z})}) is isomorphic to ${\bf
  R}$ as a topological group with respect to addition.  Observe also
that $\pi_l'$ maps the image of (\ref{({bf Z} times {bf R}) / zeta({bf
    Z})}) in $A_S / \Delta_S({\bf Z})$ onto $C_l$, since we already
know that $\pi_l$ maps the image of $Z \times {\bf R}$ in $A_S$ onto
$B_l$.

        Let $E_l$ be the subgroup of ${\bf Z} \times {\bf Z}$
generated by $(p^l \, {\bf Z}) \times \{0\}$ and $\zeta({\bf Z})$,
which is the same as
\begin{equation}
\label{E_l = {(x, y) in {bf Z} times {bf Z} : x - y in p^l {bf Z}}}
        E_l = \{(x, y) \in {\bf Z} \times {\bf Z} : x - y \in p^l \, {\bf Z}\}.
\end{equation}
It is easy to see that
\begin{equation}
\label{({bf Z} times {bf R}) / E_l}
        ({\bf Z} \times {\bf R}) / E_l
\end{equation}
is isomorphic to ${\bf R} / p^l \, {\bf Z}$ as a topological group
with respect to addition.  Of course, ${\bf R} / p^l \, {\bf Z}$ is
isomorphic to ${\bf R} / {\bf Z}$ as a topological group, which is
isomorphic to the multiplicative group ${\bf T}$ of complex numbers
with modulus $1$.  Because $B_l$ is isomorphic as a topological ring to
\begin{equation}
\label{({bf Z} times {bf R}) / ((p^l {bf Z}) times {0})}
        ({\bf Z} \times {\bf R}) / ((p^l \, {\bf Z}) \times \{0\}),
\end{equation}
$C_l$ is isomorphic as a topological group to (\ref{({bf Z} times {bf
    R}) / E_l}), and hence also to ${\bf T}$.  One can also get this
by considering the composition of $\pi_l'$ with the natural embedding
of (\ref{({bf Z} times {bf R}) / zeta({bf Z})}) into $A_S /
\Delta_S({\bf Z})$, which is a continuous homomorphism from (\ref{({bf
    Z} times {bf R}) / zeta({bf Z})}) onto $C_l$ whose kernel is equal
to the image of $(p^l \, {\bf Z}) \times \{0\}$ in (\ref{({bf Z} times
  {bf R}) / zeta({bf Z})}) under the natural quotient mapping from
${\bf Z} \times {\bf R}$ onto (\ref{({bf Z} times {bf R}) / zeta({bf
    Z})}).

        Suppose now that $\phi$ is any continuous homomorphism from
$A_S / \Delta_S({\bf Z})$ into ${\bf T}$.  Let $U$ be the set of
$z \in {\bf T}$ with positive real part, which is a relatively open
set in ${\bf T}$ that contains $1$.  Thus $\phi^{-1}(U)$ is an open
set in $A_S / \Delta_S({\bf Z})$ that contains $0$.  This implies that
$\phi$ maps
\begin{equation}
\label{q_S'((p_1^{l_1} Z_{p_1}) ... (p_k^{l_k} Z_{p_k}) times {0})}
        q_S'((p_1^{l_1} \, {\bf Z}_{p_1}) \times \cdots \times
                                    (p_k^{l_k} \, {\bf Z}_{p_k}) \times \{0\})
\end{equation}
into $U$ for some $l$, where $q_S'$ is the natural quotient mapping
from $A_S$ onto $A_S / \Delta_S({\bf Z})$, as usual.  It follows that
(\ref{q_S'((p_1^{l_1} Z_{p_1}) ... (p_k^{l_k} Z_{p_k}) times {0})})
is contained in the kernel of $\phi$ for some $l$, since the only
subgroup of ${\bf T}$ contained in $U$ is the trivial group $\{1\}$.
This permits $\phi$ to be expressed as the composition of $\pi_l'$
with a continuous homomorphism from $C_l$ into ${\bf T}$.  Of course,
the latter can be characterized in the usual way, because $C_l$ is
isomorphic to ${\bf T}$ as a topological group.

\end{document}